\documentclass{agtart_a}
\pdfoutput=1
\usepackage{pinlabel}
\title {Invariants of curves in $\mathbb{R}P^2$ and $\mathbb{R}^2$}

\author{Abigail Thompson}
\givenname{Abigail}
\surname{Thompson}
\address{Mathematics Department\\
         University of California\\\newline
         Davis, CA 95616\\
         USA}
\email{thompson@math.ucdavis.edu}
\urladdr{}

\volumenumber{6}
\issuenumber{}
\publicationyear{2006}
\papernumber{76}
\startpage{2175}
\endpage{2187}

\doi{}
\MR{}
\Zbl{}


\received{2 February 2006}
\revised{}
\accepted{2 May 2006}
\published{19 November 2006}
\publishedonline{19 November 2006}
\proposed{}
\seconded{}
\corresponding{}
\editor{CPR}
\version{}

\arxivreference{math.GT/0602003}

\AtBeginDocument{\let\bar\wbar\let\tilde\wtilde\let\hat\what}

\makeatletter
\def\cnewtheorem#1[#2]#3{\newtheorem{#1}{#3}[section]
\expandafter\let\csname c@#1\endcsname\c@theorem}
\makeatother

\newtheorem{theorem}{Theorem}
\cnewtheorem{prop}[theorem]{Proposition}
\cnewtheorem{lemma}[theorem]{Lemma}
\cnewtheorem{claim}[theorem]{Claim}
\cnewtheorem{cor}[theorem]{Corollary}
\theoremstyle{definition}
\cnewtheorem{defin}[theorem]{Definition}
\cnewtheorem{defins}[theorem]{Definitions}
\cnewtheorem{example}[theorem]{Example}
\cnewtheorem{xca}[theorem]{Exercise}

\begin{document}

\begin{asciiabstract}
There is an elegant relation found by Fabricius-Bjerre [Math. Scand 40 (1977)
20--24] among the double tangent lines, crossings, inflections points,
and cusps of a singular curve in the plane.  We give a new
generalization to singular curves in RP^2.  We note that the
quantities in the formula are naturally dual to each other in RP^2,
and we give a new dual formula.
\end{asciiabstract}

\begin{webabstract} 
There is an elegant relation found by Fabricius-Bjerre [Math. Scand 40 (1977)
20--24] among the double tangent lines, crossings, inflections points,
and cusps of a singular curve in the plane.  We give a new
generalization to singular curves in $\mathbb{R}P^2$.  We note that the
quantities in the formula are naturally dual to each other in $\mathbb{R}P^2$,
and we give a new dual formula.
\end{webabstract}

\begin{htmlabstract}
There is an elegant relation found by Fabricius-Bjerre [Math. Scand 40 (1977)
20&ndash;24] among the double tangent lines, crossings, inflections points,
and cusps of a singular curve in the plane.  We give a new
generalization to singular curves in <b>R</b>P<sup>2</sup>.  We note that the
quantities in the formula are naturally dual to each other in <b>R</b>P<sup>2</sup>,
and we give a new dual formula.
\end{htmlabstract}

\begin{abstract}
There is an elegant relation \cite{FB2} among
the double tangent lines, crossings, inflections points,
and cusps of a singular curve in the plane.
We give a new generalization to singular curves in $\mathbb{R}P^2$.  We note that the 
quantities in the formula are naturally dual to each other in $\mathbb{R}P^2$, 
and we give a new dual formula.   
\end{abstract}
\maketitle

\section{Introduction}\label{sec1}

Let $K$ be a smooth immersed curve in the plane.  Fabricius-Bjerre 
\cite{FB1} found the following relation among the double tangent 
lines, crossings, and inflections points for a generic $K$:
$$
T_1-T_2=C+(1/2)I
$$
where $T_1$ and $T_2$ are the number of exterior and interior double 
tangent lines of $K$, $C$ is the number of crossings, and $I$ is the 
number of inflection points.  Here ``generic'' means roughly that the
interesting attributes of the curve are invariant under small smooth perturbations.    Fabricius-Bjerre remarks on an example 
due to Juel which shows that the theorem cannot be straightforwardly 
generalized to the projective plane.  A series of papers followed. 
Halpern \cite{H} re-proved the theorem and obtained some additional 
formulas using analytic techniques.   Banchoff 
\cite{B} proved an analogue of the theorem for piecewise linear 
planar curves, using deformation methods.   Fabricius-Bjerre  
gave a variant of the theorem for curves with cusps \cite {FB2}. 
Weiner \cite{W} generalized the formula to closed curves lying on a 
2--sphere.  Finally Pignoni \cite{P} generalized the formula to curves 
in real projective space, but his formula depends, both in the 
statement and in the proof, on the selection of a base point for the 
curve.  Ferrand \cite{F} relates the Fabricius-Bjerre and Weiner formulas
to Arnold's invariants for plane curves.  Note that any formula for curves
in $\mathbb{R}P^2$ is more general than one for curves in $\mathbb{R}^2$, since one
can specialize to curves in $\mathbb{R}^2$ by considering curves lying inside a
small disk in $\mathbb{R}P^2$.     

There are two main results in this paper.  The first is a 
generalization of the theorem in \cite {FB2} to $\mathbb{R}P^2$, with no 
reference to a basepoint on the curve.  The original theorem is transparently a 
special case of this result, which is not surprising as the 
techniques used to prove it are a combination of those found in \cite{FB1} and 
in \cite{W}.  The difficulties encountered in the generalization 
are due to the problems in distinguishing between two ``sides'' of a
closed geodesic in $\mathbb{R}P^2$.  These are overcome by a careful attention to the natural 
metric on the space inherited from the round 2--sphere of radius one.   

\begin{figure}[ht!]
\centering
\includegraphics[scale=0.50]{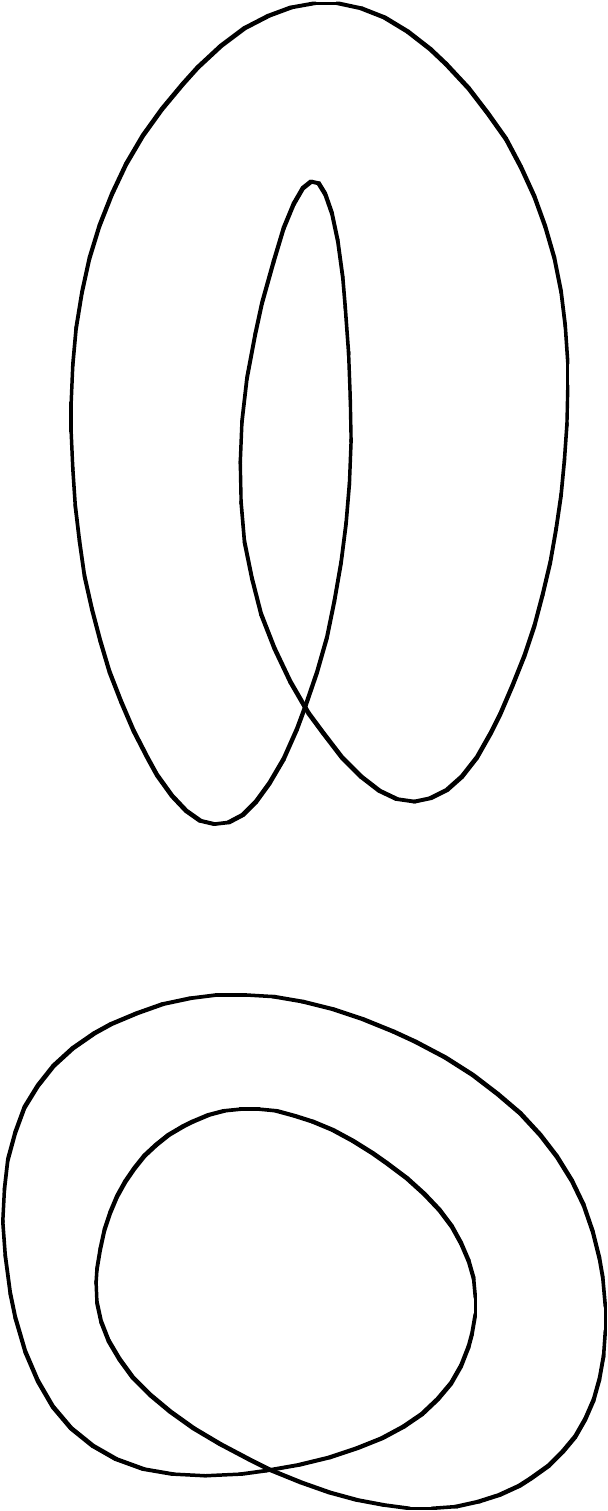}\label{twocurves}
\caption{}
\end{figure}

The main results are tied together by the observation that, in the 
version of the original formula which includes cusps \cite {FB2}, the quantities in the formula 
are naturally dual to each other in $\mathbb{R}P^2$.  This leads to the second, more surprising,  
main result, which is a dual formula for generic curves in $\mathbb{R}P^2$. 
This specializes to a new formula for generic smooth curves in the plane. 
This new formula has the interesting property that it reveals delicate 
geometric distinctions between topologically similar planar 
curves, for example quantifying some of the differences between the 
two curves shown in \fullref{twocurves}.

The outline of the paper is as follows: in \fullref{sec2} we state and 
prove the generalization of \cite{FB2} to curves in $\mathbb{R}P^2$.  In 
\fullref{sec3} we describe the duality between terms of the formula.  In 
\fullref{sec4} we state and prove the dual formulation, and give its 
corollaries for planar curves.
\vspace{-2pt}

\section{A Fabricius-Bjerre formula for curves in $\mathbb{R}P^2$}\label{sec2}
\vspace{-2pt}

Let $\mathbb{R}P^2$ be endowed with the spherical metric, inherited from its double cover, the 
round 2--sphere of radius one.  With this metric, a simple closed 
geodesic (or projective line)  in $\mathbb{R}P^2$ has length $\pi$.  The figures will use a 
standard disk model for $\mathbb{R}P^2$, in which the boundary of the disk 
twice covers a closed geodesic.
Let $K$ be a generic oriented closed curve in $\mathbb{R}P^2$, which is smoothly immersed 
except for cusps of type 1, that is, cusps at which locally the two 
branches of $K$ coming into the cusp are on opposite sides of the 
tangent geodesic.  We postpone the definition of {\em generic} until the end of section 3. 
We will need some definitions.

{\bf Definitions}\par
Let {\em $\tau_p$} be the geodesic tangent to $K$ at $p$, with 
orientation induced by $K$.

Let {\em $a_p$}, the {\em antipodal point to $p$}, be the point on 
$\tau_p$ a distance $\pi/2$ from $p$.

$\tau_p$ is divided by $p$ and $a_p$ into two pieces.
Let {\em ${\tau_p}^+$} be the segment from $p$ to $a_p$ and  {\em 
${\tau_p}^-$} the segment from $a_p$ to $p$.
At cusp points ${\tau_p}^+$ and  ${\tau_p}^-$ are not well-defined.

Let {\em $\nu_p$} be the normal geodesic to $K$ at $p$.

Let {\em $c_p$} (which lies on $\nu_p$) be the center of curvature of 
$K$ at $p$, that is, the center of the osculating circle to $K$ at $p$.

We orient {\em $\nu_p$} so that the length of the (oriented) segment 
from $p$ to $c_p$ is less than the length of the segment from $c_p$ 
to $p$.  This orientation is well-defined except at cusps and 
inflection points.

There is a natural duality from $\mathbb{R}P^2$ to itself.  Under this duality 
simple closed geodesics, or projective lines,  in $\mathbb{R}P^2$ are sent to points and vice versa. 
This duality is most easily described by passing to the 2--sphere $S$ 
which is the double cover of $\mathbb{R}P^2$; in this view a simple closed 
geodesic in $\mathbb{R}P^2$ lifts to a great circle on $S$.  If this great 
circle is called the equator,  the dual point in $\mathbb{R}P^2$ is the image 
of the north (or south) pole.

Under this duality the image of $K$ is a {\em dual curve} $K'$.   To 
describe $K'$ we need only observe that a point on $K$ comes equipped 
with a tangent geodesic, $\tau_p$.  The {\em dual point to $p$}, 
called $p'$,  is the point dual to the tangent geodesic $\tau_p$.

Another useful description is that $p'$ is the point a distance 
$\pi/2$ along the normal  geodesic to $K$ at $p$.  Notice that
$\nu_p=\nu_{p'}$ and $c_{p}=c_{p'}$.

An ordered pair of points $(p,q)$ on $K$ is an {\em antipodal pair} if $q=a_p$.

Let {\em $Y_p$} be the geodesic dual to the point $c_p$.

Let $(p,q)$ be an antipodal pair.  Then
 $Y_p$ and $\tau_p$
intersect at $q$ and divide $\mathbb{R}P^2$ into two regions, $R_1$ and $R_2$.
One of the regions, say $R_1$, contains $c_p$.  The geodesic $\tau_q$
lies in one of the two regions.     An antipodal pair $(p,q)$ is of
type 1 if  $\tau_q$ lies in $R_1$,  type 2 if  $\tau_q$ lies in
$R_2$.  Let $A_1$ be the number of type 1 antipodal pairs of $K$,
$A_2$ the number of type 2.

\begin{figure}[ht!]\small
\centering
\includegraphics[scale=0.60]{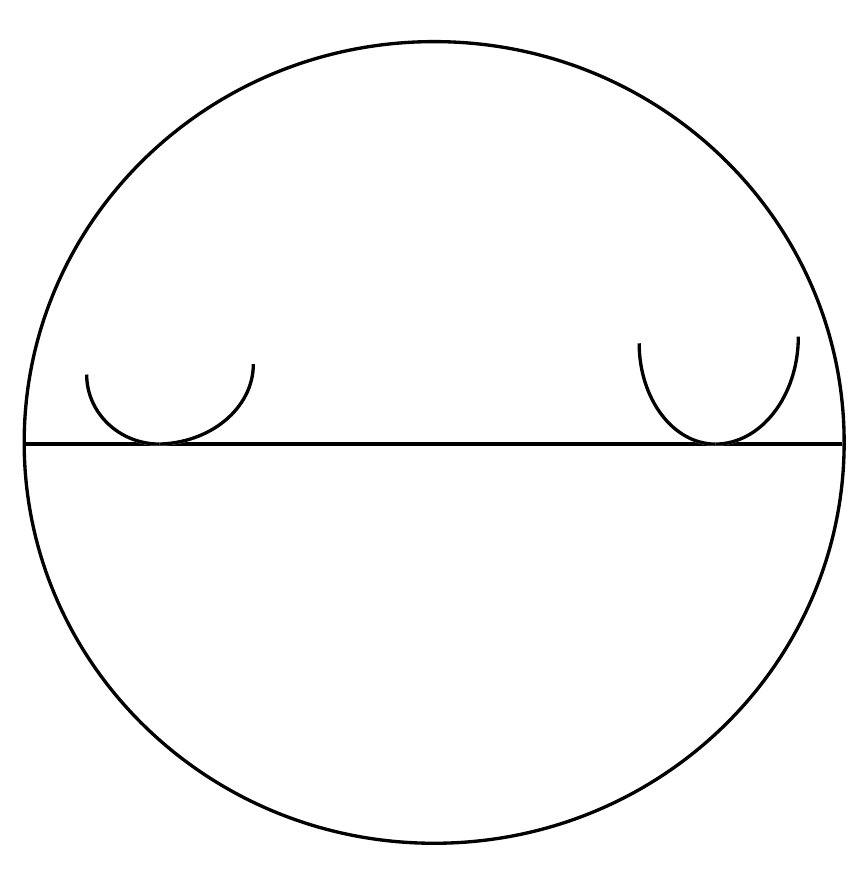}\qquad\qquad
\includegraphics[scale=0.60]{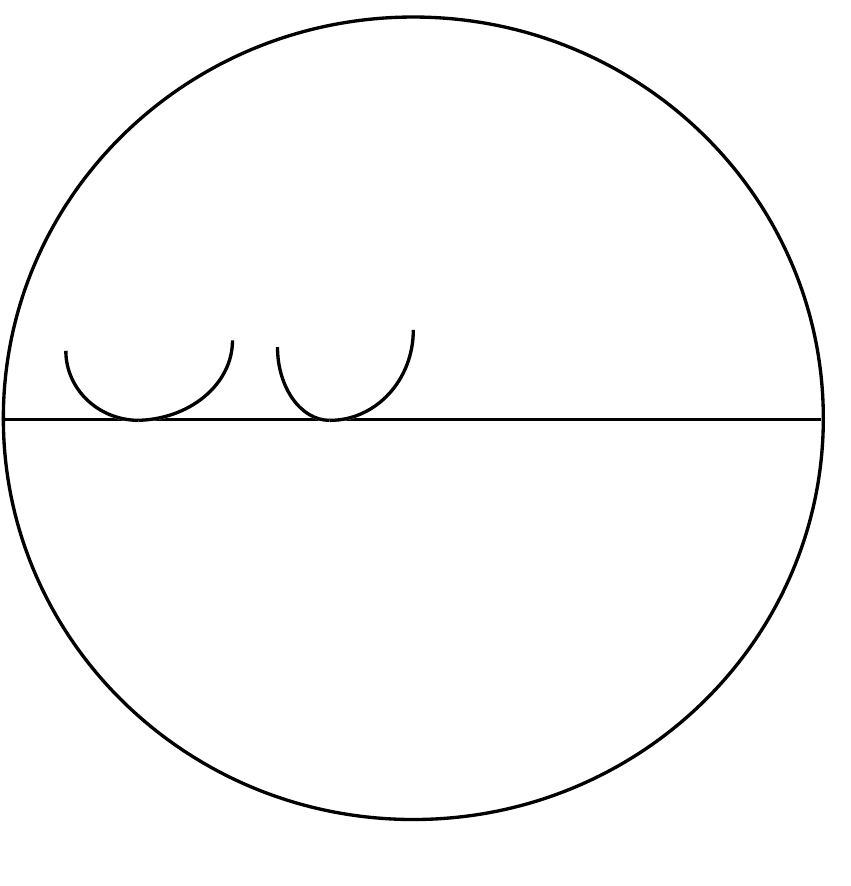}
\cl{Type 1 double supporting geodesic \qquad\qquad
Type 2 double supporting geodesic}
\caption{}\label{double tangent types}
\label{fig2}
\end{figure}

$T$ is a {\em double-supporting geodesic} of $K$ if $T$  is either a 
double tangent geodesic, a tangent geodesic through a cusp or a 
geodesic through two cusps.  The two tangent or cusp points of $K$ 
divide $T$ into two segments, one of which has length less than 
$\pi/2$.  We distinguish two types of double supporting geodesics, 
depending on whether the two points of $K$ lie on the same side of 
this segment  (type 1)  or opposite sides (type 2).  Let
{\em $T_1$} be the number of double supporting geodesics of $K$ of 
type 1, {\em $T_2$} the number of type 2 (see \fullref{fig2}).

\begin{figure}[ht!]\small
\centering
\includegraphics[scale=0.70]{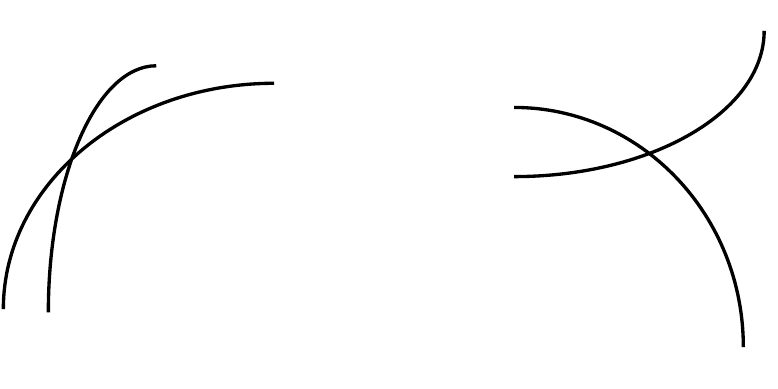}\label{crossing types}
\cl{type 1 crossing\qquad\qquad type 2 crossing}
\caption{}
\label{fig3}
\end{figure}

The tangent geodesics at a crossing of $K$ define four angles, two of 
which, $\alpha$ and $\beta$, are less than $\pi/2$.  In a small 
neighborhood of a crossing there are four segments of $K$.   The 
crossing is of type 1 if one of these segments lies in $\alpha$ and 
another in $\beta$, type 2 if two lie in $\alpha$ or two lie in 
$\beta$.  Let $C_1$ be the number of type 1 crossings of $K$, $C_2$ 
the number of type 2 (see \fullref{fig3}). 
\vspace{2pt}

\vspace{2pt}

Let {\em $I$} be the number of inflection points of $K$.
\vspace{2pt}

Let {\em $U$} be the number of (type 1) cusps of $K$.
\vspace{2pt}

We are now ready to state the first main theorem, which is a 
generalization of the main theorem of \cite {FB1} to the projective 
plane.  We note that (unlike \cite{P}) we do not need to 
choose a base-point for $K$.
\vspace{2pt}

\begin{theorem}\label{main1}
Let $K$ be a generic singular curve in $\mathbb{R}P^2$ with type 1 cusps.  Then
$$
T_1-T_2=C_1+C_2+(1/2)I+U-(1/2)A_1+(1/2)A_2
$$
\end{theorem}

\begin{proof}
The proof proceeds as in \cite {FB2}, with some caution being 
required at antipodal pairs and at cusp points.  We choose a starting 
point $p$ on $K$. Let ${M_p}^+$ be the number of times $K$ intersects 
${\tau_p}^+$,  ${M_p}^-$ be the number of times $K$ intersects 
${\tau_p}^-$, and $M_p={M_p}^+- {M_p}^-$.  We keep track of how $M_p$ 
changes as we traverse the knot.   Double-supporting geodesics, 
crossings, cusps and inflection points all behave as in \cite{FB2}. 
Suppose $p$ is a point of an antipodal pair $(p,q)$. Let $p_1$ be a 
point immediately before $p$ on $K$, $p_2$ a point immediately after. 
If $(p,q)$ is of type 1 then $\tau_{p_1}$ intersects the arc of $K$ 
containing $q$ on ${\tau_ {p_1}}^-$ and  $\tau_ {p_2}$ intersects the 
arc of $K$ containing $q$ on ${\tau_ {p_2}}^+$, hence $M_p$ increases 
by 2 as we pass through $p$. If $(p,q)$ is of type 2 then $\tau_ 
{p_1}$ intersects the arc of $K$ containing $q$ on ${\tau_ {p_1}}^+$ 
and  $\tau_ {p_2}$ intersects the arc of $K$ containing $q$ on 
${\tau_{p_2}}^-$, hence $M_p$ decreases by 2 as we pass through $p$. 
This is easiest to see by approximating $K$ near $p$ by a circle 
centered at $c_p$.  As we traverse a piece of this circle from $p_1$ 
through $p$ to $p_2$, the antipodal points to $p_1$ and $p_2$ lie on 
the geodesic $ Y_p$.  Hence the critical distinction to be made at $q$ 
is where the tangent geodesic at $q$ lies in relation to $Y_p$ and 
$\tau_p$.   This is the exactly the distinction between type 1 and 
type 2 antipodal pairs.
\end{proof}

\section{Duality in $\mathbb{R}P^2$}\label{sec3}

We describe the dual relations between crossings and double 
tangencies, cusps and inflection points, and antipodal points and 
normal-tangent pairs (defined below).

\begin{figure}[ht!]
\centering
\labellist\small
\pinlabel $q$ [t] <0pt, -2pt> at 92 332
\pinlabel $c(p)$ [t] at 155 332
\pinlabel $p$ [tl] <0pt, -2pt> at 212 332
\pinlabel {$(p,q)$ is a type 1 normal-tangent pair} at 146 254
\pinlabel $q$ [t] <0pt, -2pt> at 185 156
\pinlabel $c(p)$ [t] at 155 156
\pinlabel $p$ [tl] <0pt, -2pt> at 212 156
\pinlabel {$(p,q)$ is a type 2 normal-tangent pair} at 146 87
\endlabellist
\includegraphics[scale=0.70]{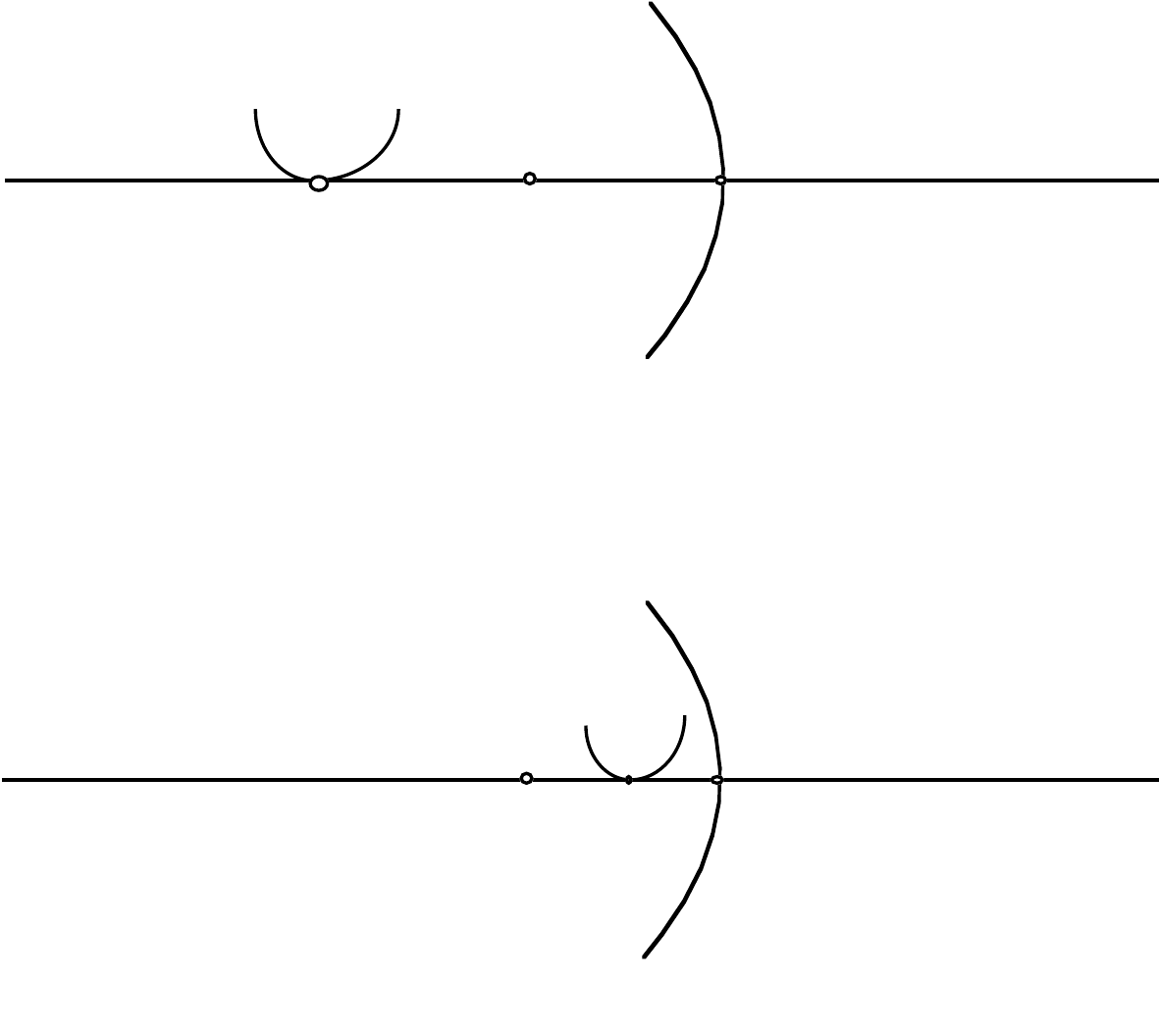}\label{normal-tangent}
\caption{}
\label{fig4}
\end{figure}

\medskip
{\bf Definitions}\qua
The points $p$ and $c_p$ divide $\nu_p$ into two pieces, {\em 
${\nu_p}^+$} from $p$ to $c_p$ and  {\em ${\nu_p}^-$} from $c_p$ to 
$p$.  An ordered pair of points $(p,q)$ on $K$ is a {\em 
normal-tangent pair} if $\tau_q=\nu_p$.
A normal-tangent pair $(p,q)$ is of type 1 if $q$ lies on ${\nu_p}^-$, 
type 2 if $q$ lies on ${\nu_p}^+$ (\fullref{fig4}).  Let $N_1$ be the number of type 
1 normal-tangent pairs of $K$, $N_2$ the number of type 2 .

\begin{prop} \label{dualcor}
Let $K$ be a generic curve in $\mathbb{R}P^2$, with dual curve $K'$.  Let $i=1,2$.
Then:
\begin{enumerate}
\item A crossing of type $i$ in $K$ is dual to a double supporting 
geodesic of type $i$ in $K'$.
\item A cusp in $K$ is dual to an inflection point in $K'$.
\item An antipodal pair of type $i$ in $K$ is dual to a normal-tangent 
pair  of type $i$ in $K'$.
\end{enumerate}

As the dual of $K'$ is again $K$, these correspondences work in both 
directions.
\end{prop}

\begin{proof}
The proof is by construction in $\mathbb{R}P^2$.
\end{proof}

This correspondence breaks down slightly when we consider double 
supporting geod\-esics between cusps and tangents, or cusps and cusps. 
Fabricius-Bjerre suggests a small local alteration of $K$ to 
understand his argument at a cusp point, replacing the cusp point by 
a small ``bump".  Just as the dual to a small round circle in $\mathbb{R}P^2$ 
is a (long) curve that is close to a geodesic, his local change at 
cusps induces a more global change at inflection points, and so in 
order to incorporate curves with inflection points we need to add 
{\em inflection geodesics} to our picture of $K$.

\begin{figure}[ht!]\small
\centering
\includegraphics[scale=0.70]{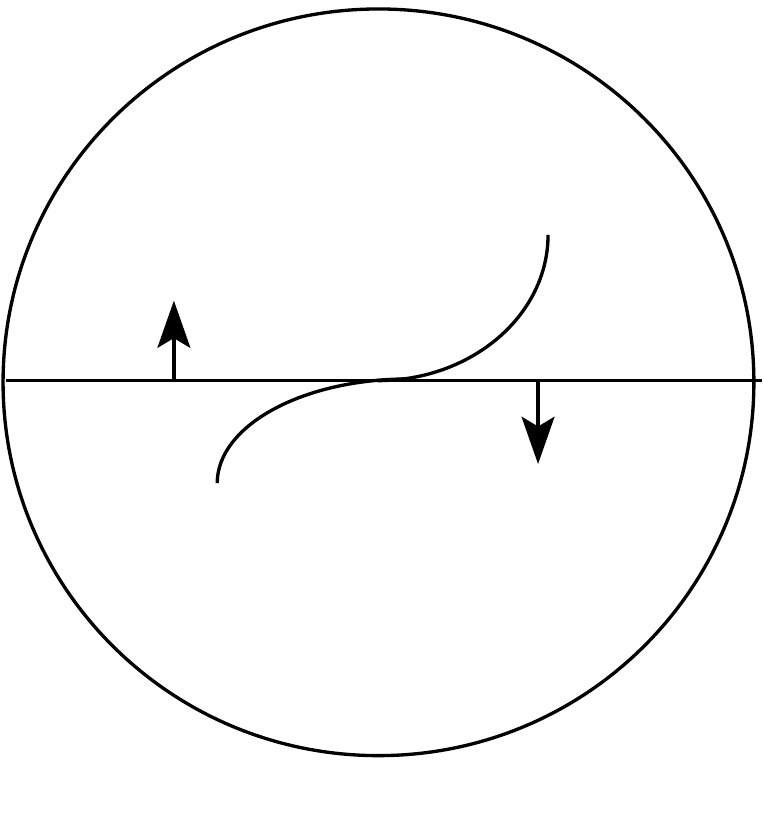}\label{inflectiongeo}
\cl{inflection geodesic with normal direction}
\caption{}
\label{fig5}
\end{figure}

Let $p$ be an inflection point of $K$, with tangent geodesic 
$\tau_p$.  Endow $\tau_p-p$ with a normal direction at each point 
(except the inflection point) by the convention shown in \fullref{fig5}.

\medskip
{\bf Definition}\qua
Call this the {\em inflection geodesic to $K$ at $p$}.

For crossings between $K$ and an inflection geodesic $\tau_p$, or 
between two inflection geodesics, the piece of $\tau_p$ in the 
neighborhood of the crossing should be construed as bending slightly 
towards its normal direction for the purposes of classifying the 
crossing type.   This convention preserves the correct duality between crossing type in $K$ and 
double supporting geodesic type in $K'$. A point on the inflection geodesic has center of 
curvature a distance $\pi/2$ in the normal direction, at $p'$. For 
$\alpha$ a point on an inflection geodesic $\tau_p$, $\nu_{\alpha}$ 
is the geodesic through $\alpha$ and $p'$.

\medskip
{\bf Definition}\qua
Denote by  $\wbar{K}$,  $K$ together with all its inflection 
geodesics.    Crossings and normal tangencies are counted as described 
above. The inflection points of $K$ (where the inflection geodesic
intersects the curve) will still be counted as simply inflection points in
$\wbar{K}$, not as new crossing points.

If $K$ is a generic curve with dual $K'$, then double supporting geodesics in $K$ involving cusp points correspond to crossings in
$\wbar{K}'$ involving inflection geodesics,  and an antipodal pair 
$(p,q)$ with $p$ a cusp point will correspond to a normal-tangent 
pair $(p',q')$ with $p'$ a point on an inflection geodesic in 
$\wbar{K}'$.   

We end this section with the definition of what it means for $K$ to be generic:

\medskip
{\bf Definition}\qua
$K$ is {\em generic} if:
\begin{itemize}
\item $K$ has a finite number of crossings, double tangent lines, cusps, inflection points, antipodal pairs, and normal-tangent pairs.
\item  Tangent geodesics at self-intersections of $K$ are neither 
parallel nor perpendicular.
\item The tangent geodesic through an 
inflection point or at a cusp is everywhere else transverse to $K$.
\item A geodesic goes through at most two tangent points or cusps of 
$K$.
\item No crossings occur at inflection points.
\item A geodesic normal to $K$ at one point is tangent to $K$ at at most one point and everywhere else transverse to $K$.
\item The distance between two points on a double-supporting geodesic is not $\pi/2$.
\item If $(p,q)$ is an antipodal pair, let $ Y_p$ be the geodesic dual 
to $c_p$.  Then $\tau_q$ should be neither $\tau_p$ nor $ Y_p$.
\item If $(p,q)$ on $K$ is a  normal-tangent pair, $q$ is not $c_p$.
\end{itemize}

\section{A dual formula, with applications}\label{sec4}

The simplest version of the dual theorem applies to curves with no 
inflection points.

\begin{theorem}\label{dualthm}
Let $K$ be a generic singular curve in $\mathbb{R}P^2$ with type 1 cusps and 
no inflection points.  Then
$$C_1-C_2=T_1+T_2+(1/2)U-(1/2)N_1+(1/2)N_2$$
\end{theorem}

Since inflection points are dual to cusps, we also have:
\begin{cor}
Let $K$ be the dual in $\mathbb{R}P^2$ of a smooth singular curve.  Then
\fullref{dualthm} holds for $K$.
\end{cor}

\proof[Proof of \fullref{dualthm}]
The theorem follows directly from duality on $\mathbb{R}P^2$, but 
it is illuminating to consider the dual of the proof of Theorem 1, 
as it provides a direct proof for curves in the plane.   In Theorem 1 
we count the number of intersections between the curve $K$ and 
$\tau_p$, with appropriate signs, as the curve is traversed once. 
Hence the main technical point is to understand the dual to $M_p$. 

\begin{figure}[ht!]
\centering
\labellist\small\hair 5pt
\pinlabel $p$ [tl] at 109 428
\pinlabel $p$ [tl] <0pt,-1pt> at 104 143
\endlabellist
\includegraphics[scale=0.50]{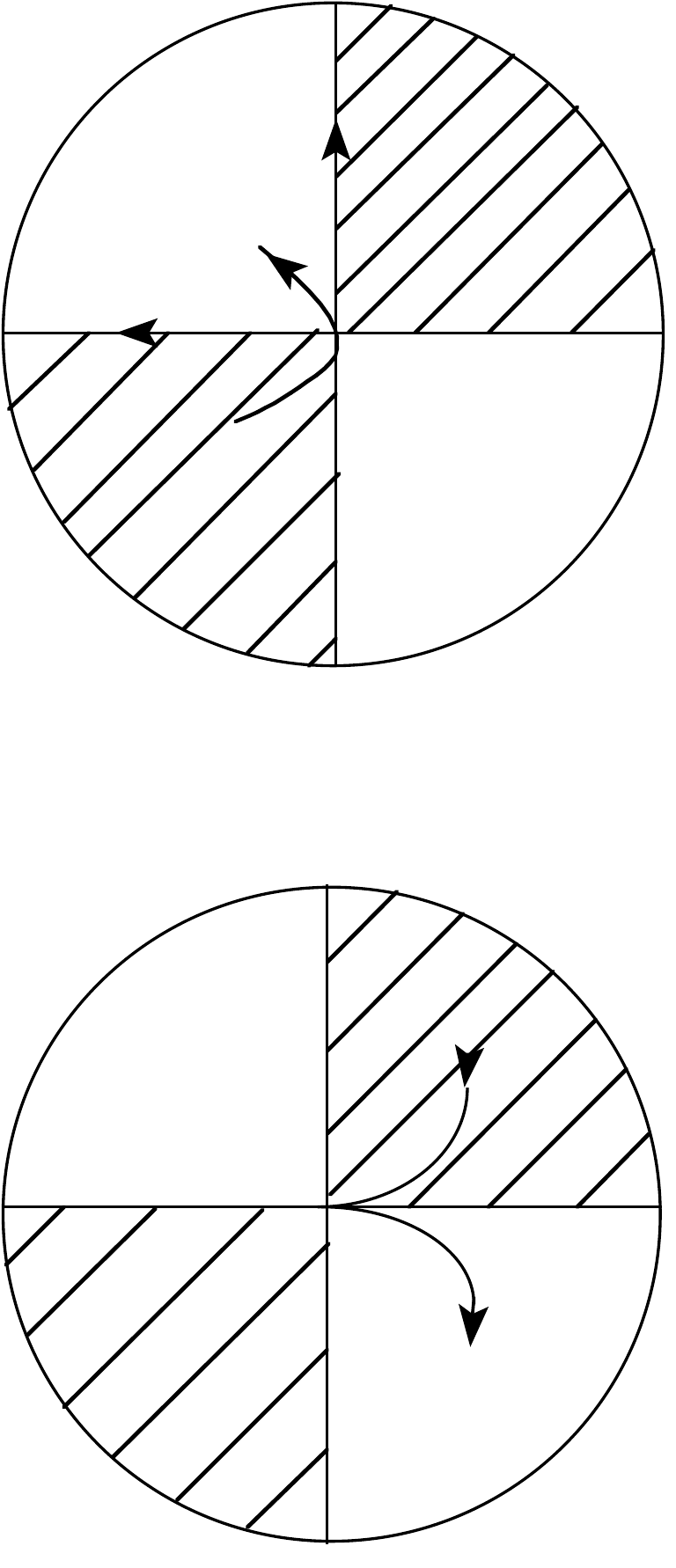}
\caption{}
\label{fig6}
\end{figure}

Assign an orientation to $K$.  The geodesics $\tau_p$ and $\nu_p$ 
intersect in a single point (at $p$) and divide $\mathbb{R}P^2$ into two 
regions.  We first define the {\em tangent-normal frame $F_p$} to $K$ 
at p as follows: $F_p$ is the union of $\tau_p$ and $\nu_p$ together 
with a black-and-white coloring of the two regions of $\mathbb{R}P^2$.  We 
color them by the convention that if we think of $\tau_p$ and $\nu_p$ 
at $P$ as being analogous to the standard $x- $ and $y-$ axes, the 
region corresponding to the quadrants where $x$ and $y$ have the same 
sign is colored white, the complementary region black (see \fullref{fig6}). 
The frame and its coloring are well-defined at points that are 
neither cusps nor inflection points.  At cusps, the orientations of 
$\tau_p$ and $\nu_p$ {\em both} reverse as we traverse $K$, with the 
happy effect that the coloring of the normal-tangent frame is 
well-defined as we pass through a cusp point (notice that this is not 
true if we allow type 2 cusps).

We now describe how the tangent-normal frame for the dual curve $K'$ is related to $M_p$  for the curve $K$.

Let $p$ be a point on $K$ with 
tangent geodesic $\tau_p$.  Let $r$ be a point of intersection 
between $K$ and $\tau_p$.  $r$ contributes either $+1$ or $-1$ to 
$M_p$, depending on where it lies relative to the antipodal point 
$a_p$.  What does $r$ correspond to in the dual picture?  Under 
duality $p$ is sent to the point $p'$, and $\nu_p=\nu_p'$.  The point 
$r$ is mapped to a geodesic $g_r$ through $p'$.  A small neighborhood 
of $r$ in $K$ is mapped to an arc of $K'$ tangent to $g_r$. If $r$ 
contributes $+1$ to $M_p$, $g_r$ lies in the white region of the 
tangent-normal frame at $p'$.  If $r$ contributes $-1$ to $M_p$, 
$g_r$ lies in the black region of the tangent-normal frame at $p'$. 
This leads us to the following definition.

\medskip
{\bf Definition}\qua
At a given point $p$ on $K$, we define $W_p$ to be the number of 
geodesics through $p$ and tangent to $K$ which lie in the white 
region and $B_p$ to be the number of geodesics through $p$ and 
tangent to $K$ which lie in the black region as defined by the 
tangent-normal frame at $p$. Let $V_p=W_p-B_p$.  The proof consists 
of tracking how $V_p$ changes as we traverse $K$ once; each type of 
singularity contributes to $V_p$ according to the following table:

\def\strutt{\vrule width 0pt depth 5pt height 12pt}
\hfill\begin{tabular}{||c|c||} \hline\hline
\strutt{\em Singularity} & {\em Contribution}\\ \hline
\strutt$C_1$ & $+4$ \\ \hline
\strutt$C_2$ & $-4$ \\ \hline
\strutt$T_i$ & $-4$ \\ \hline
\strutt$U$ & $-2$ \\ \hline
\strutt$N_1$ & $+2$ \\ \hline
\strutt$N_2$ & $-2$ \\ \hline
\end{tabular}\hfill\lower60pt\hbox{$\square$}

\medskip
We can use the natural duality directly for the general case:

\begin {theorem}
Let $K$ be a generic singular curve in $\mathbb{R}P^2$ with type 1 cusps. 
Then for $\wbar{K}$,
$$
C_1-C_2=T_1+T_2+(1/2)U+I-(1/2)N_1+(1/2)N_2
$$
\end {theorem}
\medskip

If $K$ is a curve with no cusps, inflection points, or antipodal pairs
(or for a smooth immersed curve in $\mathbb{R}^2$ with no inflection points),
then the pair of formulas:
\begin{eqnarray*}
T_1-T_2 & = & C_1+C_2\\
C_1-C_2 & = & T_1+T_2-(1/2)N_1+(1/2)N_2
\end{eqnarray*}
applies, and combining them we can obtain:
\medskip

\begin {cor}\label{tcurves}

For $K$ a curve with no cusps, inflection points, or 
antipodal pairs (or for a smooth immersed curve in $\mathbb{R}^2$ with no 
inflection points):
\begin{eqnarray*}
4T_1-4C_1=N_1-N_2\\
4T_2+4C_2=N_1-N_2
\end{eqnarray*}
\end {cor}
\medskip

\begin{figure}[ht!]\small
\centering
\includegraphics[scale=0.55]{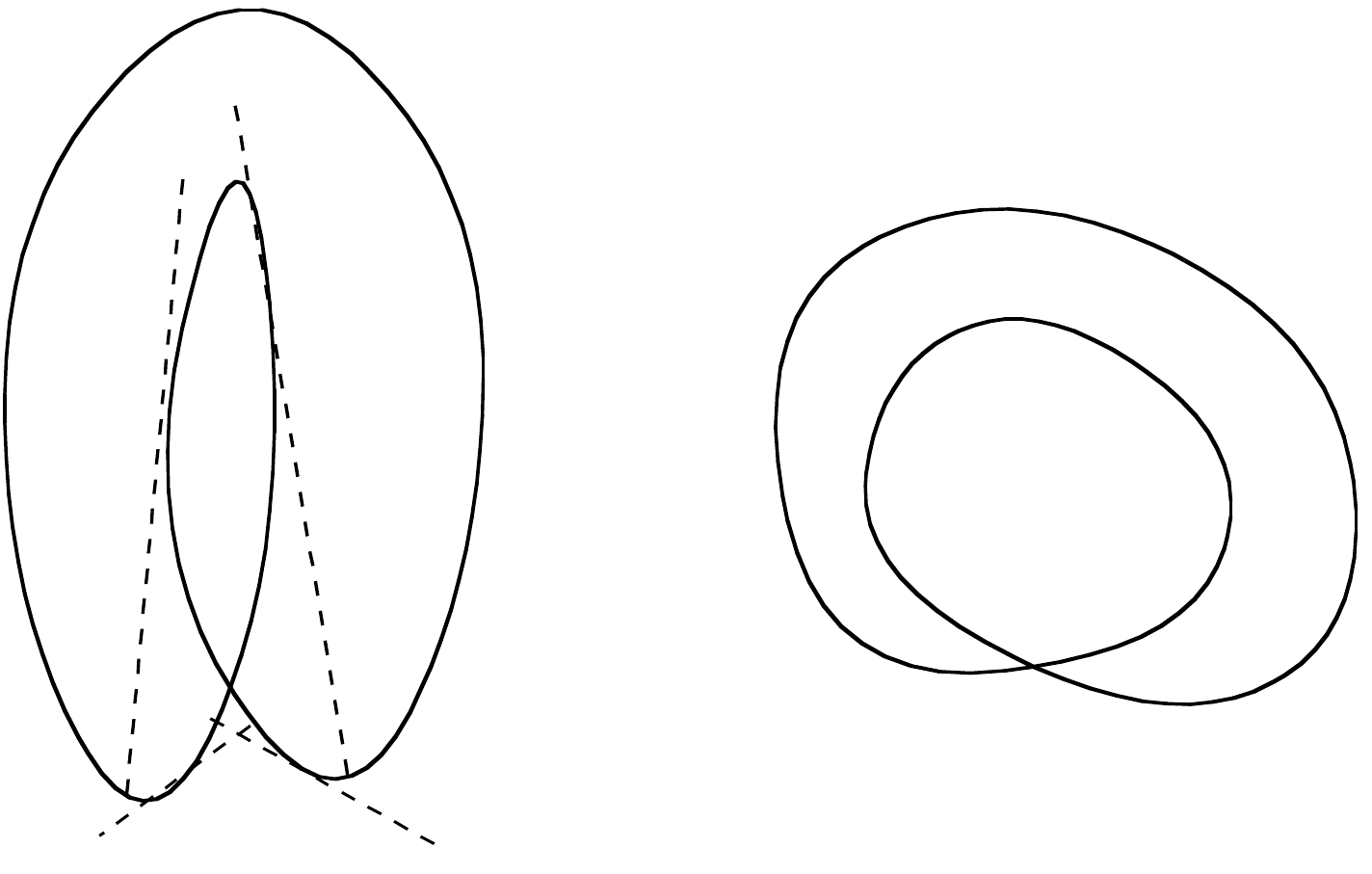}\label{twocurves*}
\cl{$C_1=0$, $C_2=1$, $N_1=4$\qquad\qquad $C_1=1$, $C_2=0$, $N_1=N_2=0$}
\caption{}\label{fig7}
\end{figure}
\medskip

Note that for the two curves shown in \fullref{twocurves} (redrawn 
in \fullref{fig7}) , we obviously have the values $T_1=1$, $T_2=0$.  For 
the right-hand curve, $C_1=1$ and $C_2=0$, while for the left, 
$C_1=0$ and $C_2=1$.   By observation, the right curve has no 
normal-tangent pairs, and the two equations in \fullref{tcurves} are easily seen to be satisfied.  Applying \fullref{tcurves} to the left-hand curve, however, we obtain
$$
4=N_1-N_2
$$
and we can locate four normal-tangent pairs of type 1 (\fullref{fig7}).
\newpage

{\bf Acknowledgement}\qua Research supported in part by an NSF grant
and by the von Neumann Fund and the Weyl Fund through the Institute
for Advanced Study.

\bibliographystyle{gtart}
\bibliography{link}

\begin{thebibliography}{}
\providecommand\bibmarginpar{\leavevmode\marginpar}
\def\urlstyle#1{{\tt #1}}

\bibitem{B}
\textbf{T\,F Banchoff},
  \href{http://links.jstor.org/sici?sici=0002-9939(197408)45:2%3C237:GGOPIT%3E%
2.0.CO%3B2-0} {\emph{Global geometry of polygons. {I}: {T}he theorem of
  {F}abricius-{B}jerre}}, Proc. Amer. Math. Soc. 45 (1974) 237--241
  \xox{MR}{0370599}

\bibitem{FB1}
\textbf{F Fabricius-Bjerre}, \emph{On the double tangents of plane closed
  curves}, Math. Scand 11 (1962) 113--116 \xox{MR}{0161231}

\bibitem{FB2}
\textbf{F Fabricius-Bjerre}, \emph{A relation between the numbers of singular
  points and singular lines of a plane closed curve}, Math. Scand. 40 (1977)
  20--24 \xox{MR}{0444673}

\bibitem{F}
\textbf{E Ferrand}, \href{http://dx.doi.org/10.1023/A:1004936711196} {\emph{On
  the {B}ennequin invariant and the geometry of wave fronts}}, Geom. Dedicata
  65 (1997) 219--245 \xox{MR}{1451976}

\bibitem{H}
\textbf{B Halpern}, \emph{Global theorems for closed plane curves}, Bull. Amer.
  Math. Soc. 76 (1970) 96--100 \xox{MR}{0262936}

\bibitem{P}
\textbf{R Pignoni}, \href{http://dx.doi.org/10.1007/BF01277967} {\emph{Integral
  relations for pointed curves in a real projective plane}}, Geom. Dedicata 45
  (1993) 263--287 \xox{MR}{1206096}

\bibitem{W}
\textbf{J\,L Weiner}, \emph{A spherical {F}abricius-{B}jerre formula with
  applications to closed space curves}, Math. Scand. 61 (1987) 286--291
  \xox{MR}{947479}

\end{thebibliography}

\end{document}